\documentclass[12pt]{amsart}
\usepackage{amssymb}
\usepackage{xypic}
\usepackage{latexsym}
\title[Jacobian of a surface with boundary]{What is the Jacobian of a
Riemann surface
with boundary?}
\author{Thomas M. Fiore and Igor Kriz}
\address{Thomas M. Fiore,
Department of Mathematics, University of Chicago, 5734 S.
University, Chicago, IL 60637, USA \\
and \\
Departament de Matem\`{a}tiques, Universitat Aut\`{o}noma de
Barcelona, 08193 Bellaterra (Barcelona), Spain}
\email{fiore@math.uchicago.edu}
\address{Igor Kriz, Department of Mathematics, University of Michigan,
2074 East
Hall, 530 Church Street, Ann Arbor, MI 48109-1043, USA}
\email{ikriz@umich.edu}
\thanks{The first authors was supported at the University of Chicago
by NSF Grant DMS-0501208. At the Universitat Aut\`{o}noma de
Barcelona he was supported by Grant SB2006-0085 of the Programa
Nacional de ayudas para la movilidad de profesores de universidad e
investigadores espa$\tilde{\text{n}}$oles y extranjeros.
The second author was supported in part by NSF Grant DMS-0305583. }

\input xypic
\newtheorem{theorem}{Theorem}[section]

\newtheorem{lemma}[theorem]{Lemma}

\theoremstyle{definition}
\newtheorem{construction}[theorem]{Construction}
\newtheorem{remark}[theorem]{Remark}
\newtheorem{definition}[theorem]{Definition}
\def\Proof{\medskip\noindent{\bf Proof: }}

\def\Z{\mathbb{Z}}

\def\C{\mathbb{C}}

\def\R{\mathbb{R}}
\def\C{\mathbb{C}}

\def\Pi{\mathbb{P}^{\infty}}

\def\qed{\hfill$\square$\medskip}

\def\Zpk{\mathbb{Z}/p^{k}}
\def\Zpk1{\mathbb{Z}/p^{k-1}}

\newcommand{\rref}[1]{(\ref{#1})}

\newcommand{\tform}[2]{\begin{array}[t]{c}#1\\
{\scriptstyle #2}\end{array}}
\newcommand{\cform}[3]{\begin{array}{c}
{\scriptstyle #3}\\
#1\\
{\scriptstyle #2}\end{array}}

\newcommand{\fracd}[2]{\frac{\displaystyle #1}{\displaystyle #2}}

\newcommand{\beg}[2]{\begin{equation}\label{#1}#2\end{equation}}
\def\r{\rightarrow}

\def\mh{\mathcal{H}}

\def\sl2{\widetilde{SL_{2}(\Z)}}

\begin{document}

\begin{abstract}
We define the Jacobian of a Riemann surface with analytically
parametrized boundary components. These Jacobians belong to a moduli
space of ``open abelian varieties'' which satisfies gluing axioms
similar to those of Riemann surfaces, and therefore allows a notion
of ``conformal field theory'' to be defined on this space. We
further prove that chiral conformal field theories corresponding to
even lattices factor through this moduli space of open abelian
varieties.
\end{abstract}

\maketitle

\vspace{5mm}

\section{Introduction}

The main purpose of the present note is to generalize the notion of
the Jacobian of a Riemann surface to Riemann surfaces with
real-analytically parametrized boundary (or, in other words,
conformal field theory worldsheets). The Jacobian of a closed
surface is an abelian variety. What structure of ``open abelian
variety'' captures the relevant data in the ``Jacobian'' of a CFT
worldsheet? If we considered Riemann surfaces with punctures instead
of parametrized boundary components, the right
answer could be easily phrased in terms of mixed Hodge structures.

\vspace{3mm}
But in worldsheets, we see more structure, and some of it is infinite-dimensional.
For example, even to a disk with analytically parametrized boundary,
one naturally assigns an infinite-dimensional symplectic form and
a restricted maximal isotropic space (cf. \cite{ps}). Any structure
we propose should certainly include such data. Additionally, in worldsheets,
boundary components can have inbound
or outbound orientation, and an inbound and outbound boundary
component can be glued to produce another worldsheet. So another test of
having the right notion of ``open abelian variety'' is that it
should enjoy a similar gluing structure.

\vspace{3mm} We should point out that it is actually a remarkably
strong requirement that a structure such as a (closed) abelian
variety could somehow be ``glued together'' from ``genus $0$'' data
similar to the situation we described above for a disk. One quickly
convinces oneself that naive approaches based on modelling somehow
the $1$-forms on a Riemann surface, together with mixed Hodge-type
integral structure data, fail to produce the required gluing. In
fact, in some sense, the desired structure must be ``pure'' rather
than ``mixed''. Note that there is no way of ``gluing'' a pure Hodge
structure out of a mixed Hodge structure which does not already
contain it: in the case of a closed Riemann surface with punctures,
the mixed Hodge structure on its first cohomogy contains the pure
Hodge structure of the original closed surface, so no gluing is
involved. Clearly, the situation is different when we are gluing a
non-zero genus surface from a genus $0$ surface with parametrized
boundary.

\vspace{3mm} There is, however, a yet stronger test. When $L$ is an
even lattice (together with a $\Z/2$-valued bilinear form $b$
satisfying a suitable condition), one has a notion of conformal
field theory associated with $L$ (\cite{scft, hk}). It could be
argued that the definition only uses additive data, so the lattice
conformal theories should ``factor through open abelian varieties''.
In some sense, if one considers the conjectured space of open
abelian varieties to be the ``Jacobian'' of the moduli space of
worldsheets (with all its structure), then one could interpret this
as a sort of ``Abelian Langlands correspondence'' for that space.
This test is also severe, as lattice conformal field theories are
known to be unexpectedly tricky. For example, the definition
of operator assigned
to a worldsheet appears to depend on the order of boundary components,
and a subtle discussion is needed to remove this (unacceptable) dependency.
This will be clarified in Section \ref{sc} below.

\vspace{3mm} In this paper, we indeed propose a notion of an open
abelian variety and answer both test questions in the affirmative.
Of course, one has to start out by being precise about what exact
abstract structure captures the notion of gluing, and then
generalize the notion of conformal field theory to be defined on
such abstract structures. Following ideas of Segal \cite{scft}, this
was done in \cite{fioret,hk,hk1}, with a correction in \cite{T}. The
desired structure is called {\em stack of pseudo commutative monoids
with cancellation} (SPCMC - see \cite{T} for a correct definition)
and a CFT is a pseudo morphism of certain SPCMC's
\beg{edcft}{\mathcal{C}\r C(\mathcal{M},\mh).} (The papers
\cite{hk,hk1} used the word ``lax'' instead of ``pseudo'', but the first
author \cite{fioret} discovered that ``pseudo''
conforms
more with existing terminology of higher category theory.)\footnote{It should
be pointed out that instead of SPCMC's, we could use
other known structures present on worldsheets which can be used for
axiomatizing CFT, for example the `cobordism approach' based on PROPs; our
structure satisfies those
axioms as well. As shown in \cite{frob}, however, when one carefully treats
the cobordism approach so no relevant axioms are omitted, the discussion
is comparable to SPCMC's.}

\vspace{3mm} In the present paper, the meaning of the target of the
map \rref{edcft}, which is defined in \cite{hk1}, plays only a
marginal role. The source of the map \rref{edcft},
however, is
important: it is the SPCMC of Segal's worldsheets. Those are
$2$-dimensional real-analytic manifolds with boundary which have a
complex structure and real-analytically parametrized boundary
components. The notion of SPCMC, which is defined in \cite{hk,hk1},
is designed to capture the operations of disjoint union and gluing
in $\mathcal{C}$, along with the fact that $\mathcal{C}$ is a
groupoid (under holomorphic maps compatible with the boundary
parametrizations), and in fact a stack over the Grothendieck
topology of complex-analytic manifolds and open covers. In
particular, gluing in $\mathcal{C}$ is defined by noticing that the
parametrized boundary components of a worldsheet can have two
possible orientations with respect to the complex structure - one
usually calls them inbound and outbound. Now from a worldsheet $X$,
another worldsheet, usually denoted by $X^\triangledown$ (despite of
the ambiguity of the symbol), can be
obtained
by gluing an inbound boundary component of $X$ to an outbound, using
the parametrizations. The notion of SPCMC is designed to capture all
the algebraic properties of these operations.

\vspace{3mm} The definition \rref{edcft} may seem mysterious, but
roughly speaking, we can imagine we have a certain finite set of
labels $A$, Hilbert spaces $H_a$ for $a\in A$, and for every
worldsheet $X$ with a map $\phi$ assigning to each boundary
component $c$ of $X$ a label $\phi(c)\in A$, a finite-dimensional
vector space $M_{X,\phi}$ and a trace class element
\beg{evac}{U_{X,\phi}\in
M_{X,\phi}\otimes\cform{\hat{\bigotimes}}{c}{}
H_{\phi(c)}^{*}\hat{\otimes} \cform{\hat{\bigotimes}}{d}{}
H_{\phi(d)} } where the tensor products are over inbound boundary
components $c$ and outbound boundary components $d$ of $X$. The
symbol $\hat{\otimes}$ means Hilbert tensor product. These elements
(called vacuum elements) are required to satisfy certain properties
which we will not list here. However, one important example is in
order. When $X$ has no boundary components (is a closed surface),
\rref{evac} becomes simply an element of $M_{X}$ ($\phi$ is dummy),
and it follows from the structure that $M_X$ is a representation of
the mapping class group $Mod(X)$.

\vspace{3mm} However, physicists noticed that in some cases (e.g.
the lattice
theories)
more is true, namely that the representation of the mapping class
group $Mod(X)$ on $M_X$ extends to the Siegel modular group
$Sp(2g,\Z)$ where $g$ is the genus of $X$ (there is a natural map
$Mod(X)\r Sp(2g,\Z)$ by taking $1$st cohomology). The question
therefore arises: what does it mean for a CFT to be
``Siegel-modular'', or, in other words, to depend only on the
cohomology of the worldsheet $X$?

\vspace{3mm} It is the main purpose of this note to provide one
possible answer to
this
question. Our approach is to define a pseudo morphism of SPCMC's
\beg{esc}{\mathcal{C}\r\mathcal{J}} where $\mathcal{J}$ is, roughly
speaking, the SPCMC of all possible `structures that look like
cohomologies of worldsheets'. We define precisely what this means,
and call such structures `open abelian varieties'.
%To give a feeling what we
%mean, recall the Siegel moduli space of abelian varieties (\cite{ab});
%this space, considered as a stack, is the ``closed sector'' of
%the SPCMC $\mathcal{S}$.

\vspace{3mm} Defining the SPCMC of open abelian varieties is our
main result. We also show that the (chiral) lattice CFT
corresponding to an even lattice indeed factors through a CFT on
$\mathcal{J}$ by the map \rref{esc}, which explains its Siegel
modularity. The reader is invited to notice that such a discussion
would be very difficult, if not impossible, if the notion of SPCMC
were not developed.

\vspace{3mm}

The present paper is organized as follows. In Section \ref{so}, we
define open abelian varieties, and discuss their moduli stack. In
Section \ref{sgl}, we discuss gluing of open abelian varieties, and
their SPCMC structure. In Section \ref{sr}, we discuss the Jacobian
map from the SPCMC of worldsheets to the SPCMC of open abelian
varieties. In Section \ref{sc}, we shall discuss the lattice
conformal field theory on the SPCMC of open abelian varieties.

\vspace{3mm}

\section{Open abelian varieties}
\label{so}

\begin{construction}
\label{const1} Let us start with the space $V_1$ of real-analytic
functions
$$f:\R\r \R$$
for which there exists a number $\Delta_f$ such that
$$f(x+2\pi)= f(x)+\Delta_f.$$
We may then alternately think of $V_1$ as a space of ``branched''
functions on $S^1$ by applying the map $e^{iz}$. There is an
antisymmetric form $S$ on $V_1$ given by
\beg{esfg}{S(f,g)=\int_{S^1}fdg-\Delta_f
g(0)-\frac{1}{2}\Delta_f\Delta_g}
(the integral over $S^1$ is interpreted as the integral from $0$ to
$2\pi$). Note that in \rref{esfg}, the term $\Delta_f g(0)$ could
have been equally well replaced by $\Delta_g f(0)$. The point is to
choose the terms so that $S(f,g)=-S(g,f)$.

\vspace{3mm} Given a pair of disjoint finite
sets $A^+$ and $A^-$, we set
$$A=A^+\amalg A^-.$$
(We think of $A^+$ as the set of
{\it outbound} and $A^-$ as the set of {\it inbound} boundary components within a connected component.) Define
$$V_A=\{f=(f_i)_i\in \cform{\prod}{i\in A^+\amalg
A^-}{}V_1|\cform{\sum}{i\in A^+}{}
\Delta_{f_i}-\cform{\sum}{i\in A^-}{}\Delta_{f_i}=0 \}/\langle
(1)_i\rangle.$$ Now choose a linear ordering on the set $A$. For
$i\in A$, define $\epsilon_i=1$ if $i\in A^+$ and $\epsilon_i=-1$ if
$i\in A^-$. Define for $f,g\in V_A$, \beg{eo1}{S_A(f,g)=
\cform{\sum}{i\in A}{}\epsilon_i S(f_i,g_i)-\frac{1}{2}
\cform{\sum}{i<j\in
A}{}\epsilon_i\epsilon_j(\Delta_{f_i}\Delta_{g_j}-
\Delta_{g_i}\Delta_{f_j}). } Note that since the space $V_A$ is
fixed, we can also give it an
integral
structure, i.e. a topological basis $B_A$ on which $S_A$ is
hyperbolic.

The exact choice does not matter. Note also that although the form
\rref{eo1} depends on the ordering of $A$, the antisymmetric forms
$S_<$ for different orderings $<$ are easily calculated from each
other, by adding differences of the corresponding terms
\beg{eodiff}{\frac{1}{2}\cform{\sum}{i<j}{}\epsilon_i
\epsilon_j(\Delta_{f_i}\Delta_{g_j}- \Delta_{g_i}\Delta_{f_j}). }
Instead of speaking of an ordering and an antisymmetric form, it
will be more useful for us to speak of a collection of antisymmetric
forms $S_<$ related to each other by the said
formulas.
We shall speak of antisymmetric forms $S_<$ related {\em in the
standard way}.
\end{construction}

\vspace{3mm}
\begin{remark}
It will be important in the sequel to note that the
form $S_<$ in fact only depends on the {\em cyclic} ordering, i.e.
if we take the smallest element $1$ of $A$ and make it the greatest
element without changing the order of the other elements, then the
form $S$ does not change. To see this, note that the operation just
described results in adding to $S$ the term
\beg{eo2}{\cform{\sum}{1\neq j\in A}{}\epsilon_1
\epsilon_j(\Delta_{f_1}\Delta_{g_j}- \Delta_{g_1}\Delta_{f_j}). }
But we are also assuming \beg{eo2a}{\cform{\sum}{j\in
A}{}\epsilon_j\Delta_{f_j}{}=0=\cform{\sum}{j\in A}{}
\epsilon_j\Delta_{g_j}, } so \rref{eo2} is equal to
$$\epsilon_{1}^{2}(\Delta_{f_1}\Delta_{g_1}-\Delta_{g_1}\Delta_{f_1})=0.$$
\end{remark}

\vspace{3mm}
\begin{remark}
\label{rem1} There is another way of relating the forms $S_{<}$,
$S_{<^\prime}$ for different orders $<$, $<^\prime$ which will be of
importance to us. Consider functions $f=(f_i)_i$ and $g=(g_i)_i$ as
above. Then define \beg{estand}{f^{\prime}_{i}=f_i-\sum\{
\epsilon_j\Delta_{f_j}|j<i \text{ and } i<^\prime j\}. } We will
refer to the map $f\mapsto f^{\prime}$ given by \rref{estand} as the
{\em standard transformation}
$$\diagram V_A\rto^\cong & V_A\enddiagram$$
corresponding to the change of the order $<$ to $<^{\prime}$. The
relation we have in mind is established by the following result.

\begin{lemma}
\label{lver} We have
$$S_{<}(f,g)=S_{<^\prime}(f^\prime,g^\prime).
$$
\end{lemma}

\Proof We have
$$
\begin{array}{l}
S_{<}(f,g)=\cform{\sum}{i\in
A}{}\epsilon_i(\int_{S^1}f_idg_i-\Delta_{f_i}g_i(0)-
\frac{1}{2}\Delta_{f_i}\Delta_{g_i})\\
-\frac{1}{2}\cform{\sum}{i<j}{}\epsilon_i\epsilon_j(\Delta_{f_i}\Delta_{g_j}
-\Delta_{g_i}\Delta_{f_j})
\\
=\cform{\sum}{i\in A}{}\epsilon_i
(\int_{S^1}f_{i}^{\prime}dg_{i}^{\prime}
+\cform{\sum}{j<i,i<^\prime
j}{}\int_{S^1}\epsilon_j\Delta_{f_j}dg_{i}^{\prime}-
\Delta_{f_{i}^{\prime}}(g_{i}^{\prime}(0)+\cform{\sum}{j<i,i<^\prime
j}{}
\epsilon_j\Delta_{g_j})\\
-\frac{1}{2}\Delta_{f^{\prime}_{i}}\Delta_{g^{\prime}_{i}})
-\frac{1}{2}\cform{\sum}{i<j}{}\epsilon_i\epsilon_j(\Delta_{f_{i}^{\prime}}
\Delta_{g_{j}^{\prime}}-\Delta_{g_{i}^{\prime}}\Delta_{f^{\prime}_{j}})
\\
=\cform{\sum}{i\in
A}{}\epsilon_i(\int_{S^1}f_{i}^{\prime}dg_{i}^{\prime}-
\Delta_{f_{i}^{\prime}}g_{i}^{\prime}(0)-\frac{1}{2}\Delta_{f_{i}^{\prime}}
\Delta_{g_{i}^{\prime}})\\
+\cform{\sum}{i\in A}{}\epsilon_i(\cform{\sum}{j<i,i<^\prime j}{}
\epsilon_j\Delta_{f^{\prime}_{j}}\Delta_{g^{\prime}_{i}})-
\cform{\sum}{i\in A}{}\epsilon_i\Delta_{f_{i}^{\prime}}
(\cform{\sum}{j<i,i<^\prime j}{}\epsilon_j\Delta_{g_{j}^{\prime}})\\
-\frac{1}{2}\cform{\sum}{i<j}{}\epsilon_i\epsilon_j(\Delta_{f^{\prime}_{i}}
\Delta_{g^{\prime}_{j}}-\Delta_{g^{\prime}_{i}}\Delta_{f^{\prime}_{j}})\\
=\cform{\sum}{i\in A}{}\epsilon_i S(f_{i}^{\prime},g_{i}^{\prime})\\
+\cform{\sum}{j<i,i<^\prime
j}{}\epsilon_i\epsilon_j(\Delta_{f_{j}^{\prime}}
\Delta_{g_{i}^{\prime}}-\Delta_{g_{j}^{\prime}}\Delta_{f_{i}^{\prime}})\\
-\frac{1}{2}\cform{\sum}{i<j}{}\epsilon_i\epsilon_j
(\Delta_{f_{i}^{\prime}}
\Delta_{g_{j}^{\prime}}-\Delta_{g_{i}^{\prime}}\Delta_{f_{j}^{\prime}})
\\
=\cform{\sum}{i\in A}{}\epsilon_i S(f_{i}^{\prime},g_{i}^{\prime})\\
-\frac{1}{2}\cform{\sum}{j<i,i<^\prime
j}{}\epsilon_i\epsilon_j(\Delta_{f_{i}^{\prime}}
\Delta_{g_{j}^{\prime}}-\Delta_{g_{i}^{\prime}}\Delta_{f_{j}^{\prime}})\\
+\frac{1}{2}\cform{\sum}{i<j,j<^\prime
i}{}\epsilon_i\epsilon_j(\Delta_{f_{i}^{\prime}}
\Delta_{g_{j}^{\prime}}-\Delta_{g_{i}^{\prime}}\Delta_{f_{j}^{\prime}})\\
-\frac{1}{2}\cform{\sum}{i<j}{} (\Delta_{f_{i}^{\prime}}
\Delta_{g_{j}^{\prime}}-\Delta_{g_{i}^{\prime}}\Delta_{f_{j}^{\prime}})\\
=\cform{\sum}{i\in A}{}\epsilon_i S(f_{i}^{\prime},g_{i}^{\prime})
-\frac{1}{2}\cform{\sum}{i<^\prime j}{}\epsilon_i\epsilon_j
(\Delta_{f_{i}^{\prime}}
\Delta_{g_{j}^{\prime}}-\Delta_{g_{i}^{\prime}}\Delta_{f_{j}^{\prime}})\\
=S_{<^\prime}(f^{\prime},g^{\prime}).
\end{array}
$$
\qed

\vspace{3mm}

It is also of interest to us that when $<^\prime$ is obtained from
$<$ by moving the greatest element $i$ to the lowest, then
$f^{\prime}_{j}=f_j$ for $j\neq i$, and $f^{\prime}_{i}$ is obtained
from $f_i$ by adding the constant function equal to
$\epsilon_i\Delta_{f_i}$. This means that when $<$ and $<^{\prime}$
correspond to the same cyclic order, the standard transformation is
not necessarily the identity, but is given by adding to each $f_i$ a
constant function which is a fixed integral multiple of
$\Delta_{f_i}$.
\end{remark}

\vspace{3mm} \noindent
\begin{definition}
\label{defo} An {\em open abelian variety} $(C,U,S,W,\iota,V^{\perp}_{\Z})$
consists of a (possibly empty) set $C$ of finite sets (called open
connected components) $A=A^+\amalg A^-$ (whose elements are called
{\it outbound} and {\it inbound boundary components} respectively), a real vector space $U$
with, for each system of linear orders $<$ of each
$A\in C$, an embedding
%\beg{eou}{{\protect \diagram U &
%\protect\cform{\bigoplus}{A\in C}{} \protect
%V_A=:V\lto_(.7){\supseteq}^(.7){\iota_<}\\
%\enddiagram}}
\beg{eou}{{\protect \diagram V:=\protect\cform{\bigoplus}{A\in C}{}
\protect
V_A\rto_(.7){\subseteq}^(.7){\iota_<} & U\\
\enddiagram}}
such that the image $\iota_<V$ is of finite codimension. (Note that
the image $\iota_<V$ does not depend on $<$.) Further, for different
choices of orders $<$ and $<^\prime$, the embeddings $\iota_<$ and
$\iota_{<^{\prime}}$ are related by composing with the standard
transformation (see Remark \ref{rem1}). Further, a nondegenerate
real symplectic form $S$ is given on $U$, and \rref{eou} maps the
form
$$\cform{\bigoplus}{A\in C}{}S_{<}\;
\text{on $\cform{\bigoplus}{A\in C}{}V_A$}$$ to $S$. (Note that by
Remark \ref{rem1}, it suffices to verify this assumption for one
$\iota_<$.)

Next, there is given a smooth (in the standard sense, see below)
complex isotropic subspace $W\subset U_{\C}$ such that
\beg{eou1}{W\oplus\overline{W}=U_{\C}, } \beg{eou2}{2i
S(x,\overline{x})>0 \;\text{for all $x\in W$} } (here $\overline{W}$
denotes the complex conjugate of $W$).

\vspace{3mm}

Additionally, there is an {\em integral structure}, which is the
following subtle data: First, there is an integral structure on the
$S$-complement $V^{\perp}$ (=annihilator) of $\iota_<V$, which means there is a
subgroup $V^{\perp}_{\Z}$ of $V^{\perp}$ on which $S$ is
isomorphic
(but not by a given isomorphism) to a hyperbolic antisymmetric form.

\vspace{3mm} Next, we impose an identification on open abelian
varieties according to the following rule. Denote by $V_{const,\Z}$
the subgroup of $V$ consisting of functions which are constant, and
have integral
value,
on every boundary component. Similarly, let $V_{deg,\Z}$ denote the
subspace of $V$ of functions which have integral degree on each
boundary component, i.e. $\Delta_{f_j}\in\Z$ for all $j\in A\in C$.
Fix a system of linear orders $<$ on each $A\in C$. Then two open
abelian varieties $(C_1,U_1,S_1,W_1,\iota_1,V^{\perp}_{\Z,1})$ and
$(C_2,U_2,S_2,W_2,\iota_2,V^{\perp}_{\Z,2})$ are identified if $C_1=C_2$, $U_1=U_2$, $S_1=S_2$, $W_1=W_2$
and the selection of the map $\iota_<$ and $V^{\perp}_{\Z}$ is subject to
the following rules:  We
require \beg{eou3}{V^{\perp}_{\Z,1}\subseteq V^{\perp}_{\Z,2}\oplus
\iota_2 V_{const,\Z},}
\beg{eou4}{(\iota_1-\iota_2)(V_{deg,\Z})\subseteq
V^{\perp}_{\Z,2}\oplus
\iota_2 V_{const,\Z}.} Note that $\iota_1-\iota_2$ is a homomorphism.
(Note that condition \rref{eou4} implies that
$\iota_1-\iota_2$ on $V$
only depends
on the degree, as it is determined by its restriction to
$V_{deg,\Z}$, and the target of that map is discrete. It then
follows that on elements of $V$ of constant degree, in particular on
$V_{const,\Z}$, $\iota_1=\iota_2$. Because of this, one can replace
$\iota_2$ by $\iota_1$, and/or $V^{\perp}_{\Z,1}$ by
$V^{\perp}_{\Z,2}$ in \rref{eou3}. Also, because of this and
\rref{eou3}, we may replace $\iota_2$ by $\iota_1$ and/or
$V^{\perp}_{\Z,2}$ by $V^{\perp}_{\Z,1}$ on the right hand side of
\rref{eou4}.)

Note also that by Remark \ref{rem1}, the choice of $<$ does not
matter in this identification, since the identification is invariant
under standard transformation. Note also that by the same remark,
fixing $<$, we may replace $\iota_<$ by $\iota_{<^\prime}$ for any
system of orders $<^{\prime}$ which defines the same cyclic order on
each $A\in C$ without changing the open abelian variety.

%, and with integral structure extending
%the integral structure on $V$ (i.e. a set
%of hyperbolic bases, extended the chosen hyperbolic bases
%on the $V_A$'s, related by integral symplectic transformations; note
%that the codimension of \rref{eou} must be an even integer $2g$,
%and since $V$ is a symplectic subspace, this said group of integral
%transformations is isomorphic to $Sp(g,\Z)$). We shall call $g$ the
%{\em genus}.

\vspace{3mm} To define smoothness of a subspace $W$, recall that we
have a standard
polarization of
$V_{\C}$ given by the isotropic subspaces $V^+$, $V^-$ of functions
on
each
copy of $S^1$ which holomorphically (resp. antiholomorphically)
extend to the unit disk (recall that polarizations do not depend on
adding or subtracting finite-dimensional subspaces). Now we mean
that the projection of $W$ to $V^{+}$ is a Fredholm operator and the
projection of $W$ to $V^-$ is a smooth operator, i.e. its singular
values (considering the Hilbert structures on $W$, $V^-$ given by
\rref{eou2} and the analogous form on $V^-$) decrease exponentially.
\end{definition}

\vspace{3mm}

Deciding which maps to call {\it morphisms} of open abelian varieties is an interesting
problem. For the purpose of the present paper, we will choose
morphisms to be only isomorphisms, which is unambiguous.

\begin{definition}
\label{dmor} An {\em isomorphism} of open abelian varieties
$$(C,U,S,W,\iota)\r
(C^{\prime},U^{\prime},S^{\prime},W^{\prime},\iota^{\prime})$$
consists of a bijection $b:C\r C^{\prime}$, and for each $A\in C$ a
bijection $b_A:A\r b(A)$ preserving inbound and outbound boundary
components, an isomorphism $\phi:U\r U^\prime$ such that
$\phi(W)=W^\prime$, $\phi$ carries $S$ to $S^\prime$, and for each
system of orders $<$ of all $A\in C$, if we denote by $<_b$ the
order induced by the system $b_A$ on $b(A)$, $b_A$ and $\phi$
conjugate $\iota_<$ to an embedding which defines the same open
abelian variety as $\iota_{<_b}^{\prime}$.

In this paper, the category of open abelian varieties will be chosen
to be the category whose objects are open abelian varieties and
whose morphisms are isomorphisms.
\end{definition}

\vspace{3mm} The identifications imposed in Definition \ref{defo}
can be viewed more systematically in the following way: Consider a
particular embedding $\iota_0:V\r U$, and a particular hyperbolic
basis of $V^{\perp}_{\Z}$.
Then we
can identify $U$ with $V\oplus V^{\perp}_{\R}$ via this embedding.
Now consider the group of all linear transformations
$$\phi:V\oplus V^{\perp}\r V\oplus V^{\perp}$$
which can be represented by $2\times 2$ matrices
$$\left(\begin{array}{cc}
\phi_{VV} & \phi_{VV^\perp}\\
\phi_{V^\perp V} & \phi_{V^\perp V^\perp}
\end{array}
\right)$$ such that the map $\phi_{V^\perp V^\perp}$ is an integral
symplectic transformation, $\phi_{V
V^{\perp}}(V^{\perp}_{\Z})\subseteq
V_{const,\Z}$,
$\phi_{V^{\perp}V}(V_{deg,\Z})\subseteq V^{\perp}_{\Z}$,
$(\phi_{VV}-Id)(V_{deg,\Z})\subseteq V_{const,\Z}$. It is easy to
check that linear transformations of this type form a discrete group,
which we denote by $Sp_{open}(V,\Z)$. This can be considered the
group of identifications of open abelian variety data.

\vspace{3mm} More precisely, let us compute the moduli space of open
abelian varieties for a given set of open connected components. From
the definition, it follows that the moduli space is
\beg{eo3}{U(W)\backslash Sp_{sm}(U)/Sp_{open}(V,\Z). } The group $U(W)$
is the Hilbert unitary group on $W$, the group $Sp_{sm}(U)$ is the
real symplectic group of $U$ which when expressed as $2\times 2$
matrices in the decomposition $W\oplus \overline{W}$, the
off-diagonal terms are smooth operators.

\vspace{3mm} It is worth noting that the group $Sp_{sm}(U)$ is in
fact also contractible, so the moduli space is a
``$K(\pi,1)$-stack''. To prove this, by Kuiper's theorem, it
suffices to show that the coset space \beg{eo4}{U(W)\backslash
Sp_{sm}(U)} is contractible. Expressing the form $S$ as a $2\times
2$ matrix as discussed above, it is of the form
$$\left(\begin{array}{cc} 0 & -iI\\iI &0\end{array}
\right),$$ so \rref{eo4} is isomorphic to the contractible space
$$\{\exp\left(\begin{array}{cc} 0 & A\\ \overline{A} &0\end{array}
\right)| \text{$A$ is symmetric smooth}\}.$$

\vspace{3mm}

\noindent
\begin{remark}
An open abelian variety with no open connected (and hence no
boundary) components is simply a real symplectic space $U$ with
integral structure and decomposition
$$U_\C=W\oplus \overline{W}$$
where $W$ is positive-definite isotropic, in other words, $S_{W\times W}=0$ and $2i
S(x,\overline{x})>0 \;\text{for all $x\in W$}$. This is equivalent data to
an abelian variety over $\C$ as in \cite{ab}.
\end{remark}

\vspace{3mm}

\section{Gluing and SPCMC structure}

\label{sgl}

\begin{theorem}
\label{t1} There exists an SPCMC structure on the set of open
abelian varieties.
\end{theorem}

\vspace{3mm} \noindent {\bf Remark:} Before embarking on this story,
let us briefly note the following curious fact: although open
abelian varieties model the notion of open connected components, it
does not model the notion of closed connected components. Moreover,
for the same reason, while one can define genus as one half of the
codimension of $V$ in $U$, the structure does not model the genus of
an individual open connected component. It is worthwhile pointing
out that one can consider a variant of our notion which would keep
track of both closed and open connected components, and would be
simply a sequence of closed and open abelian varieties with one
connected component in our sense. Such structure would also form an
SPCMC by our arguments.

\vspace{3mm} The proof of Theorem \ref{t1} will occupy the remainder
of this
section.
First, note that the stack structure over complex manifolds and
coverings follows from the moduli space remarks at the end of the
last section. Also, the operation of sum is obvious, realized by
direct sum in the obvious sense. So the main point to discuss is
gluing.

\vspace{3mm} We have the decomposition \beg{esgl1}{U\cong V\oplus
V^{\perp} } where $V^{\perp}$ is the $S$-annihilator of $V$ in $U$.
We therefore have a canonical projection given by the decomposition
\rref{esgl1} \beg{esgl2}{p:U\r V. } Composing with the projection
$q_A$ from $V$ to $V_A$ for a connected
component $A$,
we get a projection \beg{esgl3}{p_A:U\r V_A. } Composing further,
for $j\in A$, with the projection
$$q_{A,j}:V_A\r V_1/\R$$
(where $V_1$ is the space of real analytic branched functions on
$S^1$ as in Construction \ref{const1} and $\R$ is generated by the
constants), we get a projection \beg{esgl4}{p_{A,j}:U\r V_1/\R. }
All these maps of course also have complex forms, which we will
denote by the same symbol.

\vspace{3mm} Now the idea of gluing an inbound boundary component
$i\in A^-$ to an outbound boundary component $j\in B^+$, $A,B\in C$, is
to set \beg{egc1}{U^{\triangledown}=\{a\in
U|p_{A,i}(a)=p_{B,j}(a)\}/Im(V_1). } Here by $Im(V_1)$ we denote the
image of $V_1$ in $U$ by sending an element $x\in V_1$ to the sum of
$\iota_<(x_i)$ and $\iota_<(x_j)$ where $x_i$ (resp. $x_j$) is the
same function as $x$ on the $i$'th (resp. $j$'th) boundary component
and zero everywhere else. The order $<$ is selected in such
a way that
$i$ and $j$ immediately follow each other (see discussion of Cases 1
and 2 below). Then $Im(V_i)$, by our assumptions, $S$-annihilates
$\{a\in U|p_{A,i}(a)=p_{B,j}(a)\}$, so we can choose
$S^{\triangledown}$ as the form induced by $S$. However, we will
need to show that it is a non-degenerate symplectic form. To this
end, we will actually first give an independent formula for gluing
$W$, and then show that it is compatible with \rref{egc1}.

To glue $W$, we simply take \beg{egc2}{W^{\triangledown}=\{a\in
W|p_{A,i}(a)=p_{B,j}(a)\}.} Next, we will define the set of open
connected components
$C^\triangledown$
after gluing, which will give us a space $V^\triangledown$ defined
the same way as $V$, with $C$ replaced by $C^\triangledown$, and an
embedding $\iota_{<}^{\triangledown}$ after gluing corresponding to
a system of orders $<$ before gluing. Of course, one can choose the
order $<$, since for different orders the embeddings must be related
by composing with a standard transformation. Now there are two
principal cases to distinguish:

\vspace{3mm} \noindent {\bf Case 1:} $A=B$. In this case, define
$C^\triangledown$ as the set of components $E^\triangledown=E$ when
$E\neq A\in C$, and $A^\triangledown=A-\{i,j\}$ provided
$A^\triangledown\neq
\emptyset$. Next, assume $i<j<k$ for all $k\in A-\{i,j\}$. Then let
the order $<$ after gluing be given by omitting $i$, $j$. Further,
assume that the boundary component corresponding to $i$ is inbound.
We define for $x\in V^\triangledown$, $\iota_{<}^{\triangledown}(x)$
to be the projection of $\iota_<(y) $ for any $y=(y_k)\in V$ where
$y_k=x_k$ when $k\neq i,j$, and $y_i=y_j$ is arbitrary. By
definition, this embedding preserves the symplectic form.

\vspace{3mm} \noindent {\bf Case 2:} $A\neq B$. Then
$C^\triangledown$ is the set of $E^\triangledown=E$ where $A,B\neq
E\in C$, and $A^\triangledown=B^\triangledown=(A\cup B)-\{i,j\}$,
provided $A^\triangledown\neq\emptyset$. Then assume that $i$ is the
greatest element of $A$ and $j$ is the least element of $B$. Assume
again that the $i$'th boundary component is inbound. Let the
ordering on the glued connected component $(A\cup B)-\{i,j\}$ be
obtained by juxtaposing the ordering on $A-\{i\}$ before the
ordering on $B-\{j\}$. Again, for $x\in V^\triangledown$, we define
$\iota_{<}^{\triangledown}(x)$ to be the projection of $\iota_<(y) $
for any $y=(y_k)\in V$ where $y_k=x_k$ when $k\neq i,j$, and
$y_i=y_j$ is arbitrary. (Note that in this case, there is a subtlety
due to the fact that $x_k$ is only defined up to adding two
different constants for $k\in A, B$; what we mean is that the
difference of the constants is fixed by the requirement $y_i=y_j$.)
Again, we see that this embedding preserves antisymmetric forms.

\vspace{3mm} \noindent {\bf Remark:}  In the Cases 1 and 2, when
$A^\triangledown=\emptyset$, it simply gets deleted from the data
(see comments in the paragraph below Theorem \ref{t1} at the
beginning of this section). It does not affect the rest of the
gluing procedure.

\vspace{3mm} It remains to relate the formulas \rref{egc1},
\rref{egc2}, and prove that $S$ remains non-degenerate. First, since
we have complete control over the structure of $U^{\triangledown}$,
it is easy to see that
\beg{egcg}{\begin{array}{ll}g^{\triangledown}=g+1&\text{in Case 1,}\\
g^{\triangledown}=g &\text{in Case 2}\end{array}} where
$g^{\triangledown}$ denotes $1/2$ times the codimension of
$V^{\triangledown}$ in $U^{\triangledown}$. Additionally, since
$S^{\triangledown}$ is induced from $S$ (at least for the particular
choice of orderings), we know that $W^{\triangledown}\subset
U^{\triangledown}_{\C}$, $\overline{W}^{\triangledown}\subset
U^{\triangledown}_{\C}$ are isotropic and $S^{\triangledown}$-dual
to each other, so in particular
$$W^{\triangledown}\cap \overline{W}^{\triangledown}=0
$$
and thus that the natural map
\beg{egci'}{W^{\triangledown}\oplus\overline{W}^{\triangledown}\r
U_{\C}^{\triangledown} } is injective. What remains to be shown is
that, viewing \rref{egci'} as an inclusion,
\beg{egc*}{W^{\triangledown}+\overline{W}^{\triangledown}=U^{\triangledown}_{\C},
} or in other words that the map \rref{egci'} is onto. To show this,
we will take advantage of Segal's method \cite{scft} of relative
dimension. Choosing a polarization of \beg{edecomp}{V_\C=V^+\oplus
V^-} compatible with $W$ (for example as discussed in the last
section), let
$$W_0=Im (p|_W),\; \overline{W}_0=Im (p|_{\overline{W}}).$$
Denoting relative dimension with respect to the positive space $V^+$
by $dim_{V^+}$, i.e.
$$dim_{V^+}(Q)=index(\pi_Q)
$$
for $Q\subset V_\C$ where $\pi_Q:Q\r V^+$ is the projection given by
the decomposition, we get \beg{egp1}{dim (Ker(p|_W))+dim_{V^+}W_0+
dim (Ker(p|_{\overline{W}}))+dim_{V^-}\overline{W}_0=2g} (since
$W_0$ and $\overline{W}_0$ generate $V_\C$,
$$dim_{V^+}W_0+dim_{V^-}\overline{W}_0=dim(W_0\cap \overline{W}_0)).$$
But now one has \beg{egp2}{dim
(Ker(p|_{W^{\triangledown}}))+dim_{V^{\triangledown+}}
W_{0}^{\triangledown}\geq dim (Ker(p|_W))+dim_{V^+}W_0 +\epsilon }
where $\epsilon$ is $1$ in Case 1 and $0$ in Case 2 (this shift
arises because of our treatment of the constants on connected
components). Equality arises if and only if
\beg{egp3}{W^{\triangledown}_{0}+\overline{W}^{\triangledown}_{0}=
V^{\triangledown}_\C.} Similarly, we have \beg{egp2-}{dim
(Ker(p|_{\overline{W}^{\triangledown}}))+dim_{V^{\triangledown-}}
\overline{W}_{0}^{\triangledown}\geq dim
(Ker(p|_{\overline{W}}))+dim_{V^-}\overline{W}_0 +\epsilon } and
\beg{egp1w}{dim
(Ker(p|_{W^{\triangledown}}))+dim_{V^+}W_{0}^{\triangledown}+ dim
(Ker(p|_{\overline{W}^{\triangledown}}))+
dim_{V^-}\overline{W}_{0}^{\triangledown}\leq 2g+ 2\epsilon.}
Comparing \rref{egp1}, \rref{egp2}, \rref{egp2-}, \rref{egp1w}, we
see that equality must arise in \rref{egp2}, \rref{egp2-}, so we
have \rref{egp3}, which implies \rref{egc*} by \rref{egcg} and the
comment preceeding \rref{egp3}.

\vspace{3mm}

Now integral structure is discussed as follows. First of all,
$V^{\triangledown\perp}_{\Z}$ is generated by $V^{\perp}_{\Z}$ in
Case
2,
and is generated by $V^{\perp}_{\Z}$ and elements which have
integral degree on $i$ and differ by an integral value on $i,j$, and
have $0$ projection to the other boundary components (well defined
since we are in the same boundary component) in Case 1. Such
elements must generate $V^{\perp}_{\C}$ by the discussion of the
previous paragraph. Additionally, equivalence is preserved by gluing
by direct
verification.

\vspace{3mm}

To define the operations of an SPCMC as defined in \cite{T}, we need
to soup up our gluing definition to glue simultaneously several
pairs of boundary components, each consisting of one inbound and one
outbound boundary component.

Regarding the gluing of $U$ and $W$, there are obvious
generalizations of formulas \rref{egc1} and \rref{egc2} for multiple
pairs of components. The trickiest part is the discussion of the
ordering of boundary components, since in the case of multiple
boundary components, we can no longer rely on distinguishing two
cases as we did above. The
procedure
for generalizing to the case of gluing several pairs is as follows:
First, note that for an open abelian variety $\mathfrak{X}$, we can
associate
an antisymmetric form $S_{<}$ with {\em any} ordering of the entire
set of boundary components of $\mathfrak{X}$, regardless of the open
connected
components. Simply relate the forms corresponding to the orderings
in the standard way, and the embeddings $\iota_<$ by composing with
the standard transformations. (Note that even though the components
of an element in each open connected component are only defined up
to a separate additive constant, this does not affect standard
transformations.) For the operation of disjoint union,
we
simply juxtapose the order (this is possible, as permuting
cyclically the boundary components of each disjoint summand does not
change the form $S$). The general procedure for gluing is to change
the order of boundary components (while relating $S$ in the standard
way and $\iota_<$'s by composing with standard transformations) so
that all pairs of boundary components to be glued are arranged so
that the outbound component immediately follows the inbound, i.e.
the inbound is $i$'th and the outbound is $i+1$'st, if the boundary
components are indexed by integers. The key observation is that
permuting $i$ and $i+1$ past another boundary component will not
change the value of the form $S$, since the terms of \rref{eo1}
involving $i$ and $i+1$ cancel out, since $f_i$ and $f_{i+1}$ are
the same function when gluing. Similarly, the standard
transformations corresponding to such permutations are identities on
functions where $f_i$ and $f_{i+1}$ coincide. More generally,
embeddings with respect to orders of this specified form which are
related by composing with standard transformations before gluing
remain related by composing with standard transformations after
gluing, since terms coming from the glued boundary components cancel
out.

After such arrangement we take the induced embedding
$\iota^{\triangledown}_{<}$ to be associated with the order $<$
which omits all the pairs of the glued boundary components, and
leaves the order of the others unchanged. For a direct definition of
the integral structure, $V^{\triangledown\perp}_{\Z}$ is generated
by $V^{\perp}_{\Z}$ and elements which can be lifted to an element
$f$ of
the
sum of copies of $V_1$ over all the boundary components in such a
way that $f_k=0$ on any boundary component not glued, $f_i$ has
integral degree and $f_i-f_j$ is a constant integral function when
$i$, $j$ are glued. We see that this composite gluing produces an
open abelian variety, since it will be, for a particular order
selected, isomorphic to the open abelian variety obtained by gluing
the pairs of boundary components successively.

\vspace{3mm} Next, we must prove that the disjoint union and gluing
operations just defined have the coherence isomorphisms and diagrams
required in an SPCMC \cite{T}.

The coherence isomorphisms correspond simply to the identities
required for a commutative monoid with cancellation (Def 3.4 of
\cite{T}). The identities are commutativity, associativity, and
unitality of sum, unitality and transitivity of cancellation, and
distributivity of cancellation under addition. The isomorphisms are
by definition determined by what they do on $W$, where sum
corresponds to direct sum, and gluing is given by the generalization
of \rref{egc2} to multiple pairs. This is coherent with respect to
the obvious maps. It is also easy to see that the corresponding maps
are compatible with the $\iota_<$'s and the integral structure.

\vspace{3mm}

Having defined the coherence isomorphisms, we need to consider the
commutativity of coherence diagrams. Those diagrams are defined in
\cite{T}. All the diagrams required are of the following form:
Denote by $X_{a,b}$ the set of open abelian varieties with inbound
(resp. outbound) boundary components indexed by the finite set $a$
(resp. $b$). Then the basic operations are addition
$$+:X_{a,b}\times X_{c,d}\r X_{a+c,b+d}$$
and unit
$$0\in X_{0,0}$$
(here we denote the disjoint union of finite sets by $+$, and the
empty
set $0$,
which is the usual notation for commutative monoids with
cancellation),
and
gluing
$$\triangledown: X_{a+c,b+c}\r X_{a,b}.$$
We consider all {\em words} $\mathfrak{W}$ which can be written
using
$n$ distinct
variables $x_1,...,x_n$, each $x_i$ representing an open abelian
variety with inbound (resp. outbound) boundary components indexed by
$v_i$ (resp. $w_i$). The $v_i$'s and $w_i$'s are in turn words in
$m$ variables $a_1,...,a_m$ (representing finite sets), using the
finite set-level operations $+$, $0$. No variable $a_i$ is allowed
to occur more than once among the $v_i$'s, or among the $w_i$'s.
However, a variable occuring among the $v_i$'s may also occur among
the $w_i$'s (note that otherwise, the operation
$\triangledown$
could not be applied).

Now coherence diagrams \cite{T} are obtained by the following
procedure:
Alter a word $\mathfrak{W}$ repeatedly by applying one of the
identities
(commutativity, associativity, unitality of $+$, unitality and
transitivity of $\triangledown$, and distributivity). Denote the
word obtained by the end result of this sequence of alterations by
$\mathfrak{W}^\prime$. Then it is possible that the same word
$\mathfrak{W}^\prime$ could also
be obtained from
$\mathfrak{W}$ by a different sequence of alterations. Any time this
occurs, we
have an obvious corresponding coherence diagram. Our task is to show
that all such diagrams commute.

\vspace{3mm} However, this is quite easy, since an isomorphism
between open abelian
varieties
is determined by the isomorphism of the $W$'s. Now we have a
canonical
injection
\beg{einj1}{W_{\mathfrak{X}^{\triangledown}}\r W_{\mathfrak{X}},}
and
also canonical
projections \beg{einj2}{W_{\mathfrak{X}_1+\mathfrak{X}_2}\r
W_{\mathfrak{X}_i}. } Therefore, by
induction, we obtain a map \beg{einj3}{\diagram
W_{\mathfrak{W}}\protect\rto^{\phi_{\mathfrak{W}}^{i}}&
W_{\mathfrak{X}_i},
\enddiagram
} $i=1,...,n$, whose product is injective. By considering all types
of coherence isomorphisms again (units,
$\triangledown$-transitivity, $+$-commutativity and associativity),
we see that the maps \rref{einj1} and \rref{einj2} commute with the
maps induced by the coherence isomorphisms. Consequently, the two
paths $p_1$ and $p_2$ from the word $\mathfrak{W}$ to the word
$\mathfrak{W}^\prime$ induce a commutative diagram
\beg{einj4}{\diagram
W_\mathfrak{W}\ddto_{p_{j*}}\drto^{\phi_{\mathfrak{W}}^{i}}&\\
&W_i\\
W_{\mathfrak{W}^{\prime}}\urto_{\phi_{\mathfrak{W}^\prime}^{i}}, &
\enddiagram
} $j=1,2$, $i=1,...,n$. Since however the product of the maps
$\phi_{\mathfrak{W}^\prime}^{i}$ is injective, we conclude that
$p_{1*}=p_{2*}$, as required.

\vspace{3mm}

\section{The Jacobian of a worldsheet with boundary}

\label{sr}

In this paper, a {\em worldsheet} $\Sigma$ is a Riemann surface
whose boundary components $c_1,...,c_n$ are parametrized by analytic
diffeomorphisms
$$\phi_i: S^1\r c_i.$$
Taking a chart of $\Sigma$ (and thus identifying with a subset of
$\C$), boundary components oriented counterclockwise (resp.
clockwise) are called {\em inbound} (resp. {\em outbound}).
Worldsheets form an SPCMC $\mathcal{C}$, as proved in \cite{T}.

\begin{theorem}
\label{t2} There exists a morphism of SPCMC's
\beg{er1}{T:\mathcal{C}\r \mathcal{J}. } extending the Torelli map
on the moduli stack of closed Riemann surfaces.
\end{theorem}

We will also call the map $T$ the Torelli map, by extension of the
closed case. The proof of Theorem \ref{t2} will occupy the remainder
of this section.

\vspace{3mm}
\begin{definition}
A {\em cut worldsheet} is a pair $(\Sigma,\Gamma)$ where $\Sigma$ is
a worldsheet and
$$\Gamma\subset\Sigma$$
is a graph, i.e. a $1$-dimensional CW complex whose edges are
piecewise analytic, subject to the two conditions. First, the
boundary components
$c_i$ are
required to be edges of $\Gamma$ and $\phi_i(1)$ vertices (in
particular, the boundary components are not subdivided). Second, the
connected components of $\Sigma-\Gamma$ must be surfaces of genus
$0$ and their number must be equal to the number of the connected
components of $\Sigma$.
\end{definition}

\vspace{3mm} Thus, $\Gamma$ basically cuts each connected component
of $\Sigma$ into a surface of genus $0$ without disconnecting it.

\begin{lemma}
\label{lt2} A structure of a cut worldsheet (we will say simply cut
structure) exists on every worldsheet.
\end{lemma}

\Proof Without loss of generality, we can assume $\Sigma$ is
connected. To construct $\Gamma$, we can first choose a set of
disjoint collectively non-separating curves in $\Sigma$ which cut it
to a surface
$\Sigma^{\prime}$
of genus $0$, and let the vertices of $\Gamma$ be the images of $0$
under the parametrizations. Then connect the vertices by disjoint
open edges which cut $\Sigma$ into a disk. \qed

\vspace{3mm} It will be convenient to be a little more specific
about the choice of cut structure constructed in the proof of Lemma
\ref{lt2}. Note that a cut structure on a connected worldsheet
specifies a cyclic order of
boundary
components: changing for the moment the orientation of the boundary
components to outbound if necessary, this is simply the order in
which the boundary components appear if we travel the boundary of
the disk obtained by cutting the worldsheet along $\Gamma$. Now, if
$\Sigma$ is connected, we will call $(\Sigma,\Gamma)$ a {\em
standard cut structure} on $\Sigma$ if the cyclic order of the
boundary components of the genus $0$ worldsheet $\Sigma^{\prime}$
defined in the proof of Lemma \ref{lt2} is of the form
\beg{ecc}{c_1,...,c_{n},d_1,...,d_{2g},} where $c_1,...,c_n$ are the
boundary components of $\Sigma$, and $\Sigma$ is obtained from
$\Sigma^{\prime}$ by gluing $d_{2i-1}$
with
$d_{2i}$, $i=1,...,g$. We may refer to the pairs $d_{2i-1},d_{2i}$
as {\em pairs of hidden boundary components} of $\Sigma^{\prime}$. A
cut structure on a general worldsheet $\Sigma$ will be called
standard
if its restriction to every connected component of $\Sigma$ is
standard.

Now for a Riemann surface with standard cut structure
$(\Sigma,\Gamma)$, we define an open abelian variety
$T(\Sigma,\Gamma)$ as follows:

Without loss of generality, we may assume that $\Sigma$ is not
closed, for in the closed case we just take the ordinary Jacobian.
We may further assume that $\Sigma$ is connected, as there is an
obvious operation of direct sum on open abelian varieties (as
already remarked). Under the assumption, then, there is only one
open connected component
$A$,
and its elements are the boundary components of $\Sigma$. Let, then,
$W$ be
the space of holomorphic functions
$$f:\Sigma-\Gamma\r \C$$
which extend to holomorphic functions
$$\tilde{f}:\tilde{\Sigma}\r \C$$
such that for every deck transformation
$$\sigma:\tilde{\Sigma}\r\tilde{\Sigma}$$
there exists a number $n_{\sigma,f}\in\C$ such that
$$\tilde{f}(\sigma z)-\tilde{f}(z)=n_{\sigma,f}\;\text{for all
$z\in\tilde{\Sigma}$,}
$$
factored out by the space of functions constant on each connected
component.
The space $\overline{W}$ is defined analogously with the word
``holomorphic''
replaced by ``antiholomorphic''. Then we must define
$$U_{\C}=W\oplus \overline{W}.$$
To define the form $S$ on $U$, first define, for $f\in U$, a
$1$-form $\omega_f$ on $\Sigma$ by
$$\omega_f=d\tilde{f}.$$
Then define the ordering $<$ of boundary components as the order in
which
the boundary components occur on the boundary of $\Sigma-\Gamma$ in
the counterclockwise direction. (Recall that only the cyclic order
matters.)
Then define \beg{er3*}{S(f,g)=\int_{\Sigma}\omega_f\omega_g. }

\begin{lemma}
\label{ll2} The restriction \beg{er3i}{U\r V_A} is onto. More
precisely, \rref{er3i} has a splitting which is canonical on
functions of degree $0$ on each boundary component, and canonical in
the general case subject to selecting a standard cut structure on
$\Sigma$.
\end{lemma}

\Proof Assume without loss of generality that $\Sigma$ is connected
and not closed. Recall that by the Dirichlet principle, for a
(single-valued) real-analytic function $\phi_0$ on $\partial
\Sigma$,
there
exists a unique harmonic function $\phi$ on $\Sigma$ such that
$$\phi|_{\partial\Sigma}=\phi_0.$$
We can then represent uniquely
$$\phi\in W\oplus \overline{W},$$
which gives a canonical splitting of \rref{er3i} on functions of
degree
$0$.
To find a splitting on functions of non-zero degrees, note that,
using the notation \rref{ecc}, $c_1$,..., $c_n$, $d_{2}$,
$d_4$,...,$d_{2g}$
and the paths $p_1,...,p_g$ on the boundary of the disk $D$ from the
vertex $v_i$ of $\Gamma$ on $d_{2i-1}$ and the corresponding point
on $d_{2i}$ form a basis of $H_1(\Sigma,\Z)$. Therefore, there
exists a harmonic form with any
given
residues along $c_{1},...,c_{n}$ with sum $0$, and residues $0$
along $d_2,...,d_{2g}$, $p_1,...,p_g$. Integrating the form we
obtain a
function
$\phi$, and subtracting $\phi$ from the original function reduces
the general case to the degree $0$ case in terms of existence and
uniqueness.
\qed

\begin{lemma}
\label{ll1} Let $(\Sigma, \Gamma)$ be a genus $0$ cut worldsheet and
let $<$ be an order of boundary components compatible with the
cyclic order specified by the cut. Then, for real analytic functions
$f,g$ on $\partial \Sigma$,
\beg{ell1}{S(\tilde{f},\tilde{g})=S_A(f,g) } where $S_A$ is the form
defined in Construction \ref{const1}, $S$ is \rref{er3*}, and
$\tilde{f}, \tilde{g}$ are the harmonic continuations of $f,g$ to
the disk obtained by cutting $\Sigma$ along $\Gamma$.
\end{lemma}

\Proof Let $D$ be the disk obtained from $\Sigma$ by cutting along
$\Gamma$. By Stokes' theorem, we have
\beg{eiii}{\begin{array}{l}S(\tilde{f},\tilde{g})=
\int_{D}\omega_{\tilde{f}}\omega_{\tilde{g}}= \int_{\partial
D}\tilde{f}d\tilde{g}.
\end{array}
} We claim that the right hand side is equal to \rref{eo1} in the
order specified. To see this, we can assume that all the boundary
components are outbound, and the graph $\Gamma$ has no vertices
except the
vertices
$v_1,...,v_n$ on the boundary components $c_1,...,c_n$, and edges
connecting
$v_i, v_{i+1}$, $i=1,...,n-1$ (since we can always reach such case
by continuous deformation which does not change the value of
\rref{eiii}).

In this case, denoting by $f_i$, $g_i$ the restrictions of $f,g$ to
$c_i$, the contribution to the right hand side of \rref{eiii} other
than from the boundary components $c_1,...,c_n$ is
$$\begin{array}{l}
(g_2(0)-g_1(0)-\Delta_{g_1})\Delta_{f_1} +
(g_3(0)-g_2(0)-\Delta_{g_2})(\Delta_{f_1}+\Delta_{f_2})+...\\
...+(g_n(0)-g_{n-1}(0)-\Delta_{g_{n-1}})(\Delta_{f_1}+...+\Delta_{f_{n-1}})\\
=-\protect\cform{\sum}{i=1}{n}g_i(0)\Delta_{f_i}-
\protect\cform{\sum}{i\leq j}{}\Delta_{f_i}\Delta_{g_j}\\
=-\protect\cform{\sum}{i=1}{n}g_i(0)\Delta_{f_i}-\fracd{1}{2}
\cform{\sum}{i=1}{n}\Delta_{f_i}\Delta_{g_i}- \fracd{1}{2}\protect
\cform{\sum}{i<j}{}(
\Delta_{f_i}\Delta_{g_j}-\Delta_{g_i}\Delta_{f_j}).
\end{array}
$$
\qed

\begin{lemma}
\label{ll3} The conclusion of Lemma \rref{ll1} extends to all
worldsheets with standard cut structure, provided
$$\tilde{f}|_{d_{2i-1}}=\tilde{f}|_{d_{2i}}\;\; \text{of degree
$0$}$$
and
$$\tilde{g}|_{d_{2i-1}}=\tilde{g}|_{d_{2i}}\;\; \text{of degree
$0$}.$$
\end{lemma}

\Proof It suffices to assume, without loss of generality, that
$\Sigma$ is connected. Then simply apply Lemma \rref{ll1} to
$\Sigma^{\prime}$. The additional terms related to $d_{2i-1}$,
$d_{2i}$ cancel out. \qed

Note that the function $\tilde{f}$ in Lemma \rref{ll3} is determined
uniquely by $f$ and $\Gamma$. Thus, fixing $\Gamma$, we can now
define an open abelian variety $T(\Sigma,\Gamma)$ by choosing $W$ as
above, and letting the map \rref{eou} be defined by the
correspondence $f\mapsto \tilde{f}$. Regarding the integral
structure, a function $f\in V^{\perp}$ is
integral if
all the numbers $n_{\sigma,f}$ are integers. By the proof of Lemma
\rref{ll2},
this is equivalent to putting
$$V^{\perp}_{\Z}=\{f\in U|\; f|_{\partial\Sigma}=0,\;
deg(f|_{d_{2i}})\in \Z,\; f|_{d_{2i-1}}-f|_{d_{2i}}\in\Z\}.$$

\vspace{3mm} To show correctness of our definition, it remains to
show that $T(\Sigma,\Gamma)$ does not depend on the choice of
standard cut structure $\Gamma$. In other words, we need to show
that the open abelian varieties constructed by two different choices
$\Gamma_1$, $\Gamma_2$ of $\Gamma$ are related by conditions
\rref{eou3} and \rref{eou4}. Let us use the same notation as in
\rref{eou3} and \rref{eou4}, with $\iota_i$, $V^{\perp}_{\Z,i}$
constructed from $\Gamma_i$. Assume again, without loss of
generality, that $\Sigma$ is connected. Looking first at
\rref{eou3}, we see from the above comments that for $f\in
V^{\perp}_{\Z,1}$, $df$ has integral periods with respect to
$H_1(\Sigma,\Z)$ and $f$ has $0$ degrees on the boundary components.
These conditions do not depend on $\Gamma_i$. However, there is an
additional condition that the branch of the function $f$ on the disk
$D$
obtained by cutting $\Sigma$ along $\Gamma$ has $0$ restriction to
the boundary components of $\Sigma$. We see that changing the
fundamental
domain $D$ results in possibly selecting different branches of the
function
on the boundary components of $\Sigma$, which results in adding an
integral linear combination of the periods of $df$, which are
integral
constant
functions, as claimed.

Regarding \rref{eou4}, we have already shown that the selection of
$\tilde{f}$ is canonical in case of $f$ having $0$ degrees, so we
know \rref{eou4} in this case. In the general case, again, if $f\in
V_{deg,\Z}$, then $d\iota_if$ have integral periods with respect to
generators of $H_1(\Sigma,\Z)$. In addition, the restrictions of
$f_1$ and $f_2$ to the boundary component $c_j$ differ at most by
selection of a branch (since we use different fundamental domains
for calculating the restriction), i.e. by an integral constant
function. This proves \rref{eou4}.

\vspace{3mm} To complete the proof of Theorem \ref{t2}, it remains
to show that the map $T$ is compatible with gluing. We follow again
the two cases of the definition of gluing in the previous section.

\vspace{3mm} \noindent {\bf Case 1:} $A=B$. Assume, without loss of
generality, that $\Sigma$ is connected, $\Gamma$ is a standard cut
structure on $\Sigma$, and the boundary components are
$c_1,...,c_n$, as in \rref{ecc}. Without loss of generality, then,
$\Sigma^{\triangledown}$ is obtained from $\Sigma$ by gluing
$c_{n-1}$ and $c_n$. Then the projection $\Gamma^{\triangledown}$ of
$\Gamma$ onto $\Sigma^{\triangledown}$ defines a standard cut
structure on $\Sigma^{\triangledown}$, and
\beg{etgg}{T(\Sigma,\Gamma)^\triangledown=T(\Sigma^{\triangledown},
\Gamma^{\triangledown})} by definition.

\vspace{3mm} \noindent {\bf Case 2:} $A\neq B$. Without loss of
generality, $\Sigma=\Sigma_1\amalg \Sigma_2$ and we have standard
cut structures $\Gamma_i$ on $\Sigma_i$ such that
$$\Gamma=\Gamma_1\amalg \Gamma_2, $$
and the boundary components of $\Sigma_{i}^{\prime}$ are
$$c_{i,1},...,c_{i,n_i},d_{i,1},...,d_{i,2g_i}.$$
Without loss of generality, further, we are gluing $c_{1,n_1}$ to
$c_{2,1}$. Then we obtain a standard cut structure
$\Gamma^{\triangledown}$ on $\Sigma^{\triangledown}$ by taking the
projection of $\Gamma_1\cup \Gamma_2$ and omitting the
edge
corresponding to $c_{1,n_1}$ (or equivalently, $c_{2,1}$). Again, by
definition, we then have \rref{etgg}.

\vspace{3mm} The compatibility of $T$ with disjoint union is
obvious, as is
compatibility
with coherence isomorphisms (the point here, again, being that
isomorphisms of open
abelian varieties
are determined by the isomorphisms of the $W$'s, so the more subtle
structure
does not need to be discussed to prove commutativity of diagrams).

\vspace{3mm}

\section{The lattice conformal field theory on the SPCMC of open
abelian varieties}

\label{sc}

We begin by the same considerations as in \cite{hk}, starting on p.
351. Consider an even lattice $L$ and a bilinear form
$$b:L\times L\r\Z/2$$
which satisfies
$$b(x,x)\equiv \frac{1}{2}\langle x,x\rangle \mod 2.$$
Let $T=L_\C/L$. We let $T_{S^1}$ denote the space of all real
analytic
maps
$S^1\r T$. We choose a universal cover $T_{S^1}^{\prime}$ of
$T_{S^1}$, which can be considered as a space of maps $[0,1]\r
L_{\C}$. On $T^{\prime}_{S^1}$, we have a cocycle
$$c(\tilde{f},\tilde{g})=\exp\frac{2\pi i}{2}
\oint_{S^1}(\tilde{f}d\tilde{g}-\Delta_{\tilde{f}}g(0)+b(\Delta_{\tilde{f}},
\Delta_{\tilde{g}}))$$ but $L$ is canonically a normal subgroup of
the resulting
$\C^{\times}$-central
extension $\tilde{T}^{\prime}_{S^1}$, so we obtain a canonical
$\C^{\times}$-central extension
$\tilde{T}_{S^1}=\tilde{T}^{\prime}_{S^1}/L$,
$$1\r \C^{\times}\r \tilde{T}_{S^1}\r T_{S^1}\r 1.$$
For $\lambda\in L^{\prime}/L$ where $L^{\prime}$ is the dual of $L$,
there
is now a level $1$ Hilbert representation $\mathcal{H}_{\lambda}$ of
$\tilde{T}_{S^1}$ (its real subgroup acts by unitary bounded
operators) distinguished by the fact that the constant subgroup
$T\subset \tilde{T}_{S^1}$ acts by $e^{2\pi i \langle
?,\lambda
\rangle}$. Our conformal field theory associated with $L,b$ has
$L^{\prime}/L$ as its set of labels and $\mathcal{H}_{\lambda}$ as
its Hilbert spaces.

\vspace{3mm} Now consider an open abelian variety
$Y=(C,U,S,W,\iota)$. Assume without loss of generality that there is
only one open connected component $A$. Consider the pullback
\beg{ecft1}{\xymatrix{\tilde{W} \ar[r] \ar[d] & W \ar[d] \\
\underset{j \in A}{\bigoplus} V_1 \ar[r] & V_A} } (``putting back
the constants''). Assuming there is only one connected component,
\rref{ecft1} gives a short exact sequence \beg{ecft2}{0\r
\C\r\tilde{W}\r W\r 0. } Now let $U^{0}_{\Z}\subset U_{\Z}$ be the
sum of $V^{\perp}_{\Z}$ and the lattice spanned by $1_j\in V_0\cdot
j$, $j\in A$. Then
$$W_{L}=\{w\in\tilde{W}\otimes L|S( w,u)\in L\;
\text{for every $u\in U^{0}_{\Z}$}\}/L$$ ($L\subset L_\C\subset
\tilde{W}\otimes L$ is embedded by the first map \rref{ecft2}
tensored with $L$). We note that when $Y=T(\Sigma)$ for a worldsheet
$\Sigma$, then $W_L$ is canonically identified with the space of
holomorphic functions $\Sigma\r T=L_{\C}/L$. Next, we construct a
restriction
homomorphism
\beg{ecft3}{r:W_L\r\cform{\prod}{j\in A}{}T_{S^1}. } In fact, this
map is induced simply by tensoring with $L$ the pullback
to
$\tilde{W}$ of the projection \beg{ecft4}{r^{\prime}:W\r V_{\C}. }
In fact, let us note that we can assume without loss of generality
that \beg{ecft5}{\text{\rref{ecft4} is injective.} } Otherwise, $Y$
is a direct sum of $Ker (r^{\prime})\oplus \overline{ Ker
(r^{\prime}})$ (a closed abelian variety) and its $S$-complement.

\vspace{3mm} Next, note that \beg{ecft6}{\parbox{3.5in}{The
canonical central extension $\widetilde{\underset{j\in
A}{\prod}T_{S^1}}$ canonically splits when pulled back to $W_L$.} }
But in fact, this is completely analogous to the case of surfaces
(since the data used there depend only on the Jacobian), which is
treated in \cite{hk}, formulas (58)-(61). Then in the present case,
the conformal field theory data is given by the space of fixed
points \beg{ecft7}{\left( \tform{\hat{\bigotimes}}{j\in
A}\mathcal{H}_{\lambda_j}^{(*)} \right)^{W_L} } for labels
$\lambda_j$, $j\in A$ (to simplify notation, the superscript $(*)$
stands for the dual when $j\in A^-$ and is void when $j\in A^+$).
Here $\mathcal{H}_{\lambda}$, $\lambda\in L^{\prime}/L$ are the
level $1$ irreducible representations of $\tilde{T}_{S^1}$. In the
case of a closed abelian variety $Y$, the data required are given
simply by the space of theta functions on $Y\otimes L$ (see formula
(98) of \cite{hk}).

\vspace{3mm} The main statement to prove is that the dimension of
the space \rref{ecft7} is equal to \beg{ecft8}{|L^{\prime}/L|^{g} }
where $g$ is the genus of $Y$ when we have the condition
$$\tform{\sum}{j\in A} \epsilon_j\lambda_j=0\in L^{\prime}/L$$
where $\epsilon_j=1$ resp. $-1$ when $j$ is outbound resp. inbound,
and the dimension of the space \rref{ecft7} is $0$ otherwise. To
this end, choose a ``reference'' surface $\Sigma$ of genus $0$ (i.e.
a disk in $\C$ with a collection of disjoint open disks inside it
removed) which has boundary components which match those of $Y$,
with opposite orientation. Now the beginning point is that
\beg{ecft9}{\bigoplus\{ \cform{\hat{\bigotimes}}{j\in
A}{}\mathcal{H}_{\lambda_j}|
{\cform{\sum}{j}{}\epsilon_j\lambda_j=0}{}\} } is contained in the
space of sections of the line bundle associated with the principal
bundle \beg{ecft10}{\widetilde{\tform{\prod}{j\in
A}T_{S^1}}/Hol(\Sigma,T) } over
$${\cform{\prod}{j\in A}{}T_{S^1}}/Hol(\Sigma,T)$$
(In fact, the only reason equality does not occur is convergence
issues; a proof follows from the theory of loop groups \cite{ps}, we
do not give the details.)

\vspace{3mm} So this shows that the sum of \rref{ecft7} over
$\tform{\sum}{j}\epsilon_j\lambda_j=0$ is contained in (and equal to
if we can prove a certain convergence condition) the space of
sections of the line bundle associated with the principal bundle
\beg{ecft11}{W_L\backslash\widetilde{\tform{\prod}{j\in
A}T_j}/Hol(\Sigma,T)
} over \beg{ecft12}{W_L\backslash{\cform{\prod}{j\in
A}{}T_j}/Hol(\Sigma,T).} But \rref{ecft12} is the closed abelian
variety $A$ obtained by gluing $T\Sigma$ to $Y$ tensored with $L$,
and \rref{ecft11} is the $\theta$-bundle.

\vspace{3mm} So it remains to show the convergence condition. Again,
the method is analogous to \cite{hk}, Lemma 3. One first uses the
boson-fermion correspondence to show the convergence of the ``tower
modes'' of the vacuum operator, i.e. the summand of momentum $0$.
Lemma 5 then deals with sum over different momenta. The sum over
momenta is treated exactly in the same way in the present case. To
discuss the tower modes, there is also boson-fermion correspondence
in the category of open abelian varieties. It suffices to discuss
the genus $0$ case, where on the fermionic side, the vacuum is
represented simply by the space $W$ (or more precisely its image in
the appropriate Grassmanian). But that element is smooth because we
are working in the smooth moduli space.

\end{document}